\documentclass[12pt]{article}
\usepackage{latexsym}
\usepackage{amsmath}
\usepackage{amssymb}
\usepackage{amsfonts}

\numberwithin{equation}{section}

\newtheorem{theorem}{Theorem}[section]
\newtheorem{lemma}{Lemma}[section]

\newtheorem{remark}{Remark}[section]
\newcommand{\eproof}{{\mbox{\ }~\hfill
\mbox{\large $\Box$} \par \vskip 10pt}}

\newcommand{\R}{{\mathbb R}}

\title{Quantitative uniqueness for second order elliptic operators with strongly singular coefficients}
\author{Ching-Lung Lin\thanks{Department of Mathematics, National Chung Cheng University,
Chia-Yi 62117, Taiwan.\newline Partially supported by the National
Science Council of Taiwan.} \qquad Gen Nakamura\thanks{Department of
Mathematics, Hokkaido University, Sapporo 060-0810, Japan.
\newline Partially supported
by Grant-in-Aid for Scientific Research (B)(2)(No. 14340038) of
Japan Society for Promotion of Science.} \qquad Jenn-Nan
Wang\thanks{Department of Mathematics, Taida Institute of
Mathematical Sciences, NCTS (Taipei), National Taiwan University,
Taipei 106, Taiwan. \newline Partially supported by the National
Science Council of Taiwan.}}

\date{}

\begin{document}
\maketitle

\begin{abstract}
In this paper we study the local behavior of a solution to second
order elliptic operators with sharp singular coefficients in lower
order terms. One of the main results is the bound on the vanishing
order of the solution, which is a quantitative estimate of the
strong unique continuation property. Our proof relies on Carleman
estimates with carefully chosen phases. A key strategy in the proof
is to derive doubling inequalities via three-sphere inequalities.
Our method can also be applied to certain elliptic systems with
similar singular coefficients.
\end{abstract}

\section{Introduction}\label{sec1}
\setcounter{equation}{0}

Assume that $\Omega$ is a connected open set containing $0$ in
$\R^n$ for $n\geq 2$. Let $P(x,D)=\sum_{j,k}a_{jk}(x)D_j D_k$ be an
elliptic differential operator in $\Omega$ such that $a_{jk}(0)$ is
a real symmetric matrix and $a_{jk}(x)$ is Lipschitz continuous in
$\Omega$, where $D_j=\partial/\partial x_j$, $j=1,\cdots,n$. Note
that $a_{jk}(x)$ could be complex valued at $x\neq 0$. In this paper
we consider the following second order differential inequality:
\begin{equation}\label{1.1}
|P(x,D)u|\leq \frac{C_1}{|x|^{2}}|u|+\frac{C_2}{|x|}|\nabla
u|\quad\text{in}\quad{\Omega},
\end{equation}
where $C_2$ is sufficiently small. Before proceeding to the main
discussion, we want to point out that restrictions described above
are necessary. It is well known that the Lipschitz smoothness
requirement on $a_{ij}$ is minimal for the unique continuation to
hold \cite{pl}. Counterexamples given by Alinhac \cite{al} show that
the restriction of $a_{ij}(0)$ being real is necessary for the
strong unique continuation. On the other hand, regarding the
constant $C_2$, the strong unique continuation fails for \eqref{1.1}
if $C_2$ is not small, see \cite{alba} and \cite{wo}. Finally,
simple counterexamples also show that the singular coefficients on
the right side of \eqref{1.1} are sharp for the strong unique
continuation. Under the same assumptions, the strong unique
continuation property for \eqref{1.1} was proved by Regbaoui
\cite{Reg}. But Regbaoui did not give any quantitative estimate on
the vanishing order of $u$ satisfying \eqref{1.1}. This is our main
goal in this work. The development of qualitative unique
continuation property has a long history. We do not intend to give a
summary here. We refer to the paper \cite{kota} and references
therein for more details.

Concerning about the quantitative estimate of the uniqueness for
partial differential operators, we would like to mention several
related works. Using the frequency function, Garofala and Lin
\cite{gl1}, \cite{gl2} derived a quantitative version of the strong
unique continuation for strongly second order elliptic operators. In
\cite{gl1}, they also considered $|x|^{-2}$ potentials but without
first order terms.  In \cite{gl2}, they studied full lower order
terms with certain singular coefficients, but they are not sharp.
Also in \cite{ku}, Kukavica used the frequency function to prove the
maximal vanishing order of solutions to the strong second order
elliptic operator with essentially bounded potentials. Our method in
this paper is different from those in \cite{gl1}, \cite{gl2}, and
\cite{ku}. Our key tools are Carleman estimates. Besides of the
difference in method, the differential operator $P(x,D)$ in
\eqref{1.1} is only elliptic and the coefficients on the right hand
side of \eqref{1.1} are strongly singular. None of \cite{gl1},
\cite{gl2}, and \cite{ku} dealt with the equation as \eqref{1.1}.

On the other hand, Donnelly and Fefferman \cite{dofe} applied
Carleman's technique to derive the maximal vanishing order of the
eigenfunction with respect to the corresponding eigenvalue on a
compact smooth Riemannian manifold. Also, in \cite{lin1}, Lin
applied the Carleman estimate proved by Jerison and Kenig
\cite{jeke} to derive a quantitative estimate of the strong unique
continuation property for the Schr\"odinger equation with
$L^{n/2}_{loc}$ potential. However, the methods in \cite{dofe} and
\cite{lin1} can not be applied to \eqref{1.1} with strongly singular
coefficients. The difficulty lies in the fact that all Carleman
estimates used to treat the strong unique continuation contain only
polynomial weights, which are not "singular" enough to handle sharp
singular coefficients in the lower derivatives. In this work, we
overcome this difficulty by deriving three-sphere inequalities using
slightly singular than polynomial weights. Then we proceed to derive
doubling inequalities and the bound on the vanishing order of the
solution to \eqref{1.1} by applying three-sphere inequalities
recursively.

In this paper, for brevity, we only consider the scalar second order
elliptic operator. But our method can also be applied to the case
where $P(x,D)$ is an elliptic system as
$$P(x,D)=\text{diag}(P_1(x,D),\cdots,P_{\ell}(x,D)),$$ where
$P_j(x,D)$, $j=1,\cdots,\ell$, are second order elliptic operators
with Lipschitz coefficients and satisfy that
$P_j(0,D)=\cdots=P_{\ell}(0,D)$ with real symmetric coefficients.
All methods mentioned above do not seem to work in this general
case. Finally, we would like to mention that quantitative estimates
of the strong unique continuation are useful in studying the nodal
sets of eigenfunctions \cite{dofe}, or solutions of second order
elliptic equations \cite{hasi}, \cite{lin2}, or the inverse problem
\cite{abrv}. The main results of the paper are summarized as
follows. Assume that $B_{R_0}\subset\Omega$.
\begin{theorem}\label{thm1.1}
There exists a positive number $R_1<1$ such that if
 $\ 0<r_1<r_2<r_3\leq R_0$ and $r_1/r_3<r_2/r_3<R_1$, then
\begin{equation}\label{1.2}
\int_{|x|<r_2}|u|^2dx\leq
C\left(\int_{|x|<r_1}|u|^2dx\right)^{\tau}\left(\int_{|x|<{r_3}}|u|^2dx\right)^{1-\tau}
\end{equation}
for $u\in H^1({B}_{R_0})$ satisfying \eqref{1.1} in ${B}_{R_0}$,
where $C$ and $0<\tau<1$ depend on $r_1/r_3$, $r_2/r_3$ and
$P(x,D)$.
\end{theorem}
\begin{remark}\label{rem1.0}
From the proof, it suffices to take $R_1\le 1/4$. Moreover, the
constants $C$ and $\tau$ can be explicitly written as
$C=\max\{C_0(r_2/r_1)^n,\exp(B\beta_0)\}$ and $\tau=B/(A+B)$, where
$C_0>1$ and $\beta_0$ are constants depending on $P(x,D)$ and
\begin{eqnarray*}
&&A=A(r_1/r_3,r_2/r_3)=(\log (r_1/r_3)-1)^2-(\log (r_2/r_3))^2,\\
&&B=B(r_2/r_3)=-1-2\log (r_2/r_3).
\end{eqnarray*}
The explicit forms of these constants are important in the proof of
Theorem~\ref{thm1.2}.
\end{remark}

\begin{theorem}\label{thm1.2}
There exists a constant $C$ depending on $P(x,D)$ such that if
$u\in H^1_{loc}({\Omega})$ is a nonzero solution to \eqref{1.1}
with $C_2<C$, then we can find a constant $R_2$ depending on
$P(x,D)$ and a constant $m_1$ depending on $P(x,D)$ and
$\|u\|_{L^2(|x|<{R_2^2})}/\|u\|_{L^2(|x|<{R_2^4})}$ satisfying
\begin{equation}\label{1.3}
\limsup_{R\to 0}\frac{1}{R^{m_1}} \int_{|x|<R}|u|^2 dx>0.
\end{equation}
\end{theorem}

In view of the standard unique continuation property for
\eqref{1.1} in a connected domain containing the origin, if $u$
vanishes in a neighborhood of the origin then it vanishes
identically in $\Omega$. Theorem~\ref{thm1.2} provides an upper
bound on the vanishing order of a nontrivial solution to
\eqref{1.1}. The following doubling inequality is another
quantitative estimate of the strong unique continuation for
\eqref{1.1}.

\begin{theorem}\label{thm1.3}
Let $u\in H^1_{loc}({\Omega})$ be a nonzero solution to \eqref{1.1}.
Then there exist positive constants $R_3$ depending on $P(x,D)$, and
$C_3$ depending on $P(x,D)$, $m_1$ such that if $0<r\leq R_3$, then
\begin{equation}\label{1.4}
\int_{|x|\le{2r}}|u|^2dx\leq C_3\int_{|x|\le{r}}|u|^2dx,
\end{equation}
where $m_1$ is the constant obtained in Theorem \ref{thm1.2}.
\end{theorem}

The rest of the paper is devoted to the proofs of
Theorem~\ref{thm1.1}-\ref{thm1.3}.

\section{Proof of Theorem \ref{thm1.1}}\label{sec2}
\setcounter{equation}{0}

This section is devoted to the proof of Theorem~\ref{thm1.1}. To
begin, we recall a Carleman estimate with weight
$\varphi_{\beta}=\varphi_\beta(|x|)
=\exp(\frac{\beta}{2}(\log|x|)^2)$ derived in \cite{Reg}.

\begin{lemma}{\rm \cite[Theorem~1.2]{Reg}}\label{lem2.1}
For any $\beta >0$ large enough. Let $S$ be a small neighborhood of
$0$ and $u:S\setminus\{0\}\subset \Omega \rightarrow {\mathbb R}$,
 $u\in {C^\infty_0 ({S}\setminus\{0\})}$. Then we have
\begin{equation}\label{2.1}
\begin{array}{l}
\quad \beta^3 \int \varphi^2_\beta {|x|^{-n}|u|^2 dx}
+\beta \int \varphi^2_\beta {|x|^{-n+2}|\nabla u|^2 dx}\\
\leq \tilde C_0\int {\varphi^2_\beta |x|^{-n+4}|P(x,D)u|^2 dx},
\end{array}
\end{equation}
for some positive constant $\tilde C_0$ depending only on $P(x,D)$.
\end{lemma}

\begin{remark}\label{rem2.1}
The estimate \eqref{2.1} in Lemma \ref{lem2.1} remains valid if we
assume $u\in H^2(S\setminus\{0\})$ with compact support. This can be
easily obtained by cutting off $u$ for small $|x|$ and regularizing.
\end{remark}

We now proceed to the main part of the proof. Using regularization,
Friedrich's lemma, and ellipticity of $P(x,D)$, we can see that if
$u\in H_{loc}^1(\Omega)$ satisfies \eqref{1.1} then $u\in
H^2_{loc}(\Omega\setminus\{0\})$. To begin, we first consider the
case where $0<r_1<r_2<R<1$ and $B_R\subset\Omega$. The constant $R$
will be determined later. To use the Carleman estimate \eqref{2.1},
we need to cut-off $u$. So let $\xi(x)\in C^{\infty}_0 ({\mathbb
R}^n)$ satisfy $0\le\xi(x)\leq 1$ and
\begin{equation*}
\xi (x)=
\begin{cases}
\begin{array}{l}
0,\quad |x|\leq r_1/e,\\
1,\quad r_1/2<|x|<er_2,\\
0,\quad |x|\geq 3r_2.
\end{array}
\end{cases}
\end{equation*}
Here $e=\exp(1)$. It is easy to see that for all multiindex $\alpha$
\begin{equation}\label{gradxi}
\begin{cases}
|D^{\alpha}\xi|=O(r_1^{-|\alpha|})\ \text{for all}\ r_1/e\le |x|\le r_1/2\\
|D^{\alpha}\xi|=O(r_2^{-|\alpha|})\ \text{for all}\ er_2\le |x|\le
3r_2.
\end{cases}
\end{equation}
On the other hand, repeating the proof of Corollary~17.1.4 in
\cite{Hor3}, we can show that
\begin{equation}\label{inter}
\int_{a_1r<|x|<a_2r}||x|^{|\alpha|}D^{\alpha} u|^2dx\le
C'\int_{a_3r<|x|<a_4r}|u|^2dx,\quad|\alpha|\le 2,
\end{equation}
for all $0<a_3<a_1<a_2<a_4$ such that $B_{a_4r}\subset\Omega$, where
the constant $C'$ is independent of $r$.

Noting that the commutator $[P(x,D),\xi]$ is a first order
differential operator. Applying \eqref{2.1} to $\xi u$ and using
\eqref{1.1}, \eqref{gradxi}, \eqref{inter} implies
\begin{eqnarray*}
&& \beta^3 \int_{r_1/2<|x|<er_2} \varphi^2_\beta {|x|^{-n}|u|^2 dx}
+\beta \int_{r_1/2<|x|<er_2} \varphi^2_\beta {|x|^{-n+2}|\nabla u|^2 dx}\\
&\leq& \beta^3 \int \varphi^2_\beta {|x|^{-n}|\xi u|^2 dx}
+\beta \int \varphi^2_\beta {|x|^{-n+2}|\nabla (\xi u)|^2 dx}\\
&\leq& \tilde C_0\int {\varphi^2_\beta |x|^{-n+4}|P(x,D)(\xi u)|^2 dx}\\
\end{eqnarray*}
\begin{eqnarray}\label{2.2}
&\leq& \tilde C_0\int\varphi^2_\beta (C_1^2|x|^{-n}|\xi
u|^2+C_2^2|x|^{-n+2}|\xi\nabla u|^2)dx\notag\\
&&+\tilde C_0\int\varphi^2_\beta |x|^{-n+4}\big{|}[P(x,D),\xi]u\big{|}^2 dx\notag\\
&\leq& \tilde C_1\Big{\{}\int_{r_1/2<|x|<er_2} \varphi^2_\beta
{|x|^{-n}|u|^2 dx}
+\int_{r_1/2<|x|<er_2} \varphi^2_\beta {|x|^{-n+2}|\nabla u|^2 dx}\notag\\
&& + \int_{r_1/e<|x|<r_1/2} \varphi^2_\beta {|x|^{-n}|u|^2 dx}
+\int_{r_1/e<|x|<r_1/2} \varphi^2_\beta {|x|^{-n+2}|\nabla u|^2 dx}\notag\\
&& + \int_{er_2<|x|<3r_2} \varphi^2_\beta {|x|^{-n}|u|^2 dx}
+\int_{er_2<|x|<3r_2} \varphi^2_\beta {|x|^{-n+2}|\nabla u|^2 dx}\Big{\}}\notag\\
&\leq& \tilde{C}_2\Big{\{}\int_{r_1/2<|x|<er_2} \varphi^2_\beta
{|x|^{-n}|u|^2 dx}
+\int_{r_1/2<|x|<er_2} \varphi^2_\beta {|x|^{-n+2}|\nabla u|^2 dx}\notag\\
&&\qquad + r_1^{-n}\varphi^2_\beta(r_1/e) \int_{r_1/e<|x|<r_1/2}
(|u|^2+||x|^2\nabla u|^2)dx\notag\\
&&\qquad + r_2^{-n}\varphi^2_\beta(er_2) \int_{er_2<|x|<3r_2}
(|u|^2+||x|^2\nabla u|^2) dx\Big{\}}\notag\\
&\leq& \tilde{C}_3\Big{\{}\int_{r_1/2<|x|<er_2} \varphi^2_\beta
{|x|^{-n}|u|^2 dx}
+\int_{r_1/2<|x|<er_2} \varphi^2_\beta {|x|^{-n+2}|\nabla u|^2 dx}\notag\\
&&\qquad + r_1^{-n}\varphi^2_\beta(r_1/e) \int_{r_1/4<|x|<r_1}
|u|^2dx+r_2^{-n}\varphi^2_\beta(er_2) \int_{2r_2<|x|<4r_2}
|u|^2dx\Big{\}},\notag\\
\end{eqnarray}
where $\tilde C_1$, $\tilde C_2$, and $\tilde C_3$ are independent
of $r_1$ and $r_2$. Now letting $\beta_0\geq 1$ and
$\beta\geq\beta_0\geq 2\tilde C_3$ in \eqref{2.2}, we immediately
get that
\begin{eqnarray}\label{2.3}
&&\int_{r_1/2<|x|<er_2} \varphi^2_\beta {|x|^{-n}|u|^2 dx} +
\int_{r_1/2<|x|<er_2} \varphi^2_\beta
{|x|^{-n+2}|\nabla u|^2 dx}\notag\\
&\leq&\tilde C_4\Big{\{} r_1^{-n}\varphi^2_\beta(r_1/e)
\int_{r_1/4<|x|<r_1} |u|^2 dx + r_2^{-n}\varphi^2_\beta(er_2)
\int_{2r_2<|x|<4r_2} |u|^2 dx\Big{\}},\notag\\
\end{eqnarray}
where $\tilde C_4=1/\tilde C_3$. It follows easily from \eqref{2.3}
that
\begin{equation*}
r_2^{-n}\varphi^2_\beta(r_2) \int_{r_1/2<|x|<r_2}|u|^2 dx
\end{equation*}
\begin{eqnarray}\label{2.4}
&\leq&\int_{r_1/2<|x|<er_2} \varphi^2_\beta {|x|^{-n}|u|^2 dx}\notag\\
&\leq&\tilde C_4\Big{\{} r_1^{-n}\varphi^2_\beta(r_1/e)
\int_{r_1/4<|x|<r_1} |u|^2 dx + r_2^{-n}\varphi^2_\beta(er_2)
\int_{2r_2<|x|<4r_2} |u|^2 dx\Big{\}}.\notag\\
\end{eqnarray}
Dividing $r_2^{-n}\varphi^2_\beta(r_2)$ on the both sides of
\eqref{2.4} implies
\begin{eqnarray}\label{2.5}
&&\int_{r_1/2<|x|<r_2}|u|^2 dx\notag\\
&\leq&\tilde C_4\Big{\{}
(r_2/r_1)^{n}[\varphi^2_\beta(r_1/e)/\varphi^2_\beta(r_2)]
\int_{r_1/4<|x|<r_1} |u|^2 dx\notag\\
&&\quad + [\varphi^2_\beta(er_2)/\varphi^2_\beta(r_2)]
\int_{2r_2<|x|<4r_2} |u|^2 dx\Big{\}}\notag\\
&\leq&\tilde C_5\Big{\{}
(r_2/r_1)^{n}[\varphi^2_\beta(r_1/e)/\varphi^2_\beta(r_2)]
\int_{|x|<{r_1}} |u|^2 dx\notag\\
&&\quad +(r_2/r_1)^{n}[\varphi^2_\beta(er_2)/\varphi^2_\beta(r_2)]
\int_{|x|<{4r_2}} |u|^2 dx\Big{\}},
\end{eqnarray}
where $\tilde C_5=\max\{\tilde C_4,1\}$. With such choice of
$\tilde C_5$, we see that
$$
\tilde
C_5(r_2/r_1)^{n}[\varphi^2_\beta(r_1/e)/\varphi^2_\beta(r_2)]>1
$$
for all $0<r_1<r_2$. Adding $\int_{|x|<{r_1/2}} |u|^2 dx$ to both
sides of \eqref{2.5} and choosing $r_2\leq 1/4$, we obtain that
\begin{eqnarray}\label{2.6}
&&\int_{|x|<{r_2}}|u|^2 dx\notag\\
&\leq& 2\tilde
C_5(r_2/r_1)^{n}[\varphi^2_\beta(r_1/e)/\varphi^2_\beta(r_2)]
\int_{|x|<{r_1}} |u|^2 dx\notag\\
&&+2\tilde C_5(r_2/r_1)^{n}
[\varphi^2_\beta(er_2)/\varphi^2_\beta(r_2)] \int_{|x|<1} |u|^2 dx.
\end{eqnarray}
For simplicity, by denoting
\begin{eqnarray*}
&&A=\beta^{-1}\,\log[\varphi^2_\beta(r_1/e)/\varphi^2_\beta(r_2)]=(\log r_1-1)^2-(\log r_2)^2>0,\\
&&B=-\beta^{-1}\,\log[\varphi^2_\beta(er_2)/\varphi^2_\beta(r_2)]=-1-2\log
r_2>0,
\end{eqnarray*}
\eqref{2.6} becomes
\begin{eqnarray}\label{22.66}
&&\int_{|x|<{r_2}}|u|^2 dx\notag\\
&\leq& 2\tilde C_5(r_2/r_1)^{n}\Big{\{}\exp(A\beta)\int_{|x|<{r_1}}
|u|^2 dx+\exp(-B\beta) \int_{|x|<1} |u|^2 dx\Big{\}}.\notag\\
\end{eqnarray}

To further simplify the terms on the right hand side of
\eqref{22.66}, we consider two cases. If
$$\exp{(A\beta_0)}\int_{|x|<{r_1}} |u|^2 dx<\exp{(-B\beta_0)}\int_{|x|<{1}} |u|^2 dx,$$
then we can pick a $\beta>\beta_0$ such that
$$
\exp{(A\beta)}\int_{|x|<{r_1}} |u|^2
dx=\exp{(-B\beta)}\int_{|x|<{1}} |u|^2 dx.
$$
Using such $\beta$, we obtain from \eqref{22.66} that
\begin{eqnarray}\label{2.8}
&&\int_{|x|<{r_2}}|u|^2 dx\notag\\
&\leq& 4\tilde C_5(r_2/r_1)^{n}\exp{(A\beta)} \int_{|x|<{r_1}} |u|^2dx\notag\\
&=& 4\tilde
C_5(r_2/r_1)^{n}\left(\int_{|x|<{r_1}}|u|^2dx\right)^{\frac{B}{A+B}}\left(\int_{|x|<{1}}|u|^2dx\right)^{\frac{A}{A+B}}.
\end{eqnarray}
On the other hand, if
$$ \exp{(-B\beta_0)}\int_{|x|<{1}} |u|^2dx\leq\exp{(A\beta_0)}\int_{|x|<{r_1}} |u|^2 dx,$$
then we have
\begin{eqnarray}\label{2.9}
&&\int_{|x|<{r_2}}|u|^2 dx\notag\\
&\leq&
\left(\int_{|x|<1}|u|^2dx\right)^{\frac{B}{A+B}}\left(\int_{|x|<1}|u|^2dx\right)^{\frac{A}{A+B}}\notag\\
&\leq&
\exp{(B\beta_0)}\left(\int_{|x|<{r_1}}|u|^2dx\right)^{\frac{B}{A+B}}\left(\int_{|x|<1}|u|^2dx\right)^{\frac{A}{A+B}}.
\end{eqnarray}
Putting together \eqref{2.8}, \eqref{2.9}, and setting $\tilde
C_6=\max\{4\tilde C_5(r_2/r_1)^n,\exp{(B\beta_0)}\}$, we arrive at
\begin{equation}\label{2.99}
\int_{|x|<{r_2}}|u|^2 dx \le \tilde
C_6\left(\int_{|x|<{r_1}}|u|^2dx\right)^{\frac{B}{A+B}}\left(\int_{|x|<1}|u|^2dx\right)^{\frac{A}{A+B}}.
\end{equation}

Now for the general case, we take $R_1\le 1/4$ and consider
$0<r_1<r_2<r_3$ with $r_1/r_3<r_2/r_3\le 1/4$. By scaling, i.e.
defining $\widehat{u}(y):=u(r_3y)$ and
$\widehat{a_{ij}}(y)=a_{ij}(r_3y)$, we derive from \eqref{2.99} that
\begin{equation}\label{2.10}
\int_{|y|<{r_2/r_3}}|\widehat{u}|^2 dy \leq
C(\int_{|y|<{r_1/r_3}}|\widehat{u}|^2dy)^{\tau}(\int_{|y|<1}|\widehat{u}|^2dy)^{1-\tau},
\end{equation}
where $\tau=B/(A+B)$ with
\begin{eqnarray*}
&&A=A(r_1/r_3,r_2/r_3)=(\log (r_1/r_3)-1)^2-(\log (r_2/r_3))^2,\\
&&B=B(r_2/r_3)=-1-2\log (r_2/r_3),
\end{eqnarray*}
and $C=\max\{4\tilde C_5(r_2/r_1)^n,\exp(B\beta_0)\}$. We want to
remark that $\tilde C_5$ can be chosen independent of the scaling
factor $r_3$ provided $r_3<1$. Restoring the variable $x=r_3y$ in
\eqref{2.10} gives
$$
\int_{|x|<{r_2}}|u|^2 dx \leq
C(\int_{|x|<{r_1}}|u|^2dx)^{\tau}(\int_{|x|<{r_3}}|u|^2dx)^{1-\tau}.
$$
The proof now is complete. \eproof

\section{Proof of Theorem \ref{thm1.2} and Theorem~\ref{thm1.3}}\label{sec3}
\setcounter{equation}{0}

In this section, we prove Theorem~\ref{thm1.2} and
Theorem~\ref{thm1.3}. Without loss of generality, we assume
$P(0,D)=\Delta$ by the change of coordinates.  We begin with another
Carleman estimate derived in \cite[Lemma~2.1]{Reg}: for any $u\in
{C^\infty_0 ({\mathbb R}^n \backslash {\{0\}})}$ and for any $m\in
{\{j+\frac{1}{2},j\in{\mathbb N}\}}$ we have
\begin{equation}\label{3.1}
\sum_{|\alpha|\leq 2} \int m^{2-2|\alpha|}
|x|^{-2m+2|\alpha|-n}|D^\alpha u|^2 dx\leq
 C\int {|x|^{-2m+4-n}|\Delta u|^2 dx},
\end{equation}
where $C$ only depends on the dimension $n$.

\begin{remark}\label{rem2.2}
Using the cut-off function and regularization, estimate
\eqref{3.1} remains valid for any fixed $m$ if $u\in
H^{2}_{loc}({\mathbb R}^n \backslash {\{0\}})$ with compact
support.
\end{remark}

In view of Remark~\ref{rem2.2}, we can apply \eqref{3.1} to the
function $\chi u$ with $\chi(x) \in C^{\infty}_0 ({\mathbb
R}^n\backslash {\{0\}})$. Therefore, we define $\chi(x) \in
C^{\infty}_0 ({\mathbb R}^n\backslash {\{0\}})$ such that
$$
\chi(x)=\begin{cases}
0\quad\text{if}\quad
|x|\leq \delta/3,\\
1\quad\text{in}\quad \delta/2\leq|x|\leq(R_0+1)R_0R/4=r_4R,\\
0\quad\text{if}\quad 2r_4R\leq |x| ,
\end{cases}
$$
where $\delta\le R_0^2R/4$, $R_0>0$ is a small number which will
be chosen later and $R$ is sufficiently small satisfying $0<R\leq
R_0$. Here the number $R$ is not yet fixed and is given by
$R=(\gamma m)^{-1}$, where $\gamma>0$ is a large constant which
will be chosen later. Using the estimate \eqref{3.1} and the
equation \eqref{1.1}, we can derive that
\begin{eqnarray*}
&& \sum_ {|\alpha|\leq {2}}m^{2-2|\alpha|}\int_{\delta/2\leq|x|\leq {r_4R}} |x|^{-2m+2|\alpha|-n}|D^\alpha u|^2 dx\\
&\leq&  \sum_ {|\alpha|\leq {2}}m^{2-2|\alpha|}\int |x|^{-2m+2|\alpha|-n}|D^\alpha (\chi u)|^2 dx\\
&\leq&  C\int |x|^{-2m+4-n}|\Delta (\chi u)|^2 dx\\
&\leq& C\int_{\delta/2\leq|x|\leq {r_4R}} |x|^{-2m+4-n}|\Delta u|^2
dx+
C\int_{|x|> {r_4R}} |x|^{-2m+4-n}|\Delta(\chi u)|^2 dx\notag\\
&&+C\int_{\delta/3\leq|x|\leq \delta/2} |x|^{-2m+4-n}|\Delta(\chi u)|^2 dx\notag\\
&\leq&  \hat{C}'\int_{\delta/2\leq|x|\leq {r_4R}}
|x|^{-2m+4-n}|\Delta
u-P(x,D)u|^2 dx\notag\\
&&+\hat{C}'\int_{\delta/2\leq|x|\leq {r_4R}} |x|^{-2m+4-n}|P(x,D)u|^2 dx\notag\\
&&+C\int_{|x|> {r_4R}} |x|^{-2m+4-n}|\Delta(\chi u)|^2 dx\notag
+C\int_{\delta/3\leq|x|\leq \delta/2} |x|^{-2m+4-n}|\Delta(\chi u)|^2 dx\notag\\
\end{eqnarray*}
\begin{eqnarray}\label{3.3}
&\leq&C'\sum_ {|\alpha|= {2}}r_4^2R^2\int_{\delta/2\leq|x|\leq {r_4R}} |x|^{-2m+4-n}|D^\alpha u|^2 dx\notag\\
&&+C'C_1^2\int_{\delta/2\leq|x|\leq {r_4R}} |x|^{-2m-n}|u|^2 dx
+C'C_2^2\sum_ {|\alpha|={1}}\int_{\delta/2\leq|x|\leq {r_4R}} |x|^{-2m+2-n}|D^\alpha u|^2 dx\notag\\
&&+C\int_{|x|> {r_4R}} |x|^{-2m+4-n}|\Delta(\chi u)|^2 dx
+C\int_{\delta/3\leq|x|\leq \delta/2} |x|^{-2m+4-n}|\Delta(\chi
u)|^2 dx,\notag\\
\end{eqnarray}
where the constant $C'$ depends on $n$.

By carefully checking terms on both sides of \eqref{3.3}, we now
choose $\gamma=\sqrt{C'}$ and thus
$$R=\frac{1}{\gamma m}=\frac{1}{\sqrt{C'}m}\quad\text{and}\quad
r_4^2R^2=\frac{R_0^2(R_0+1)^2}{16m^2C'}.$$ Hence, choosing $R_0<1$
(suffices to guarantee $R_0^2(R_0+1)^2/16<1/2$), $m\ge \tilde
m_0=\tilde m_0(R_0)$, and $C_2$ sufficiently small such that
$$\frac{1}{\sqrt{C'}m}\le R_0,\quad\frac{m^2}{2}>C'C_1^2,\quad\text{and}\quad 1-C'C_2^2>\frac 12,$$ we can
remove the first three terms on the right hand side of the last
inequality in \eqref{3.3} and obtain
\begin{eqnarray}\label{3.31}
&& \sum_ {|\alpha|\leq {2}}m^{2-2|\alpha|}\int_{\delta/2<|x|< {r_4R}} |x|^{-2m+2|\alpha|-n}|D^\alpha u|^2 dx\notag\\
&\le& 2C\int_{\delta/3<|x|<\delta/2}|x|^{-2m+4-n}|\Delta (\chi u)|^2
dx\notag\\
&&+2C\int_{r_4R<|x|<2r_4R}|x|^{-2m+4-n}|\Delta (\chi u)|^2 dx.
\end{eqnarray}

In view of the definition of $\chi$, it is easy to see that for
all multiindex $\alpha$
\begin{equation}\label{chii}
\begin{cases}
|D^{\alpha}\chi|=O(\delta^{-|\alpha|})\ \text{for all}\ \delta/3<|x|<\delta/2,\\
|D^{\alpha}\chi|=O((r_4R)^{-|\alpha|})\ \text{for all}\
r_4R<|x|<2r_4R.
\end{cases}
\end{equation}
Note that $R_0^2\le r_4$ provided $R_0\le 1/15$. Therefore, using
\eqref{chii} and \eqref{inter} in \eqref{3.31}, we derive
\begin{eqnarray}\label{3.4}
&&m^2(2\delta)^{-2m-n}\int_{\delta/2<|x|\le
2\delta}|u|^2dx+m^2(R_0^2R)^{-2m-n}\int_{2\delta<
|x|\le R_0^2R}|u|^2dx\notag\\
&\le& \sum_ {|\alpha|\leq {2}}m^{2-2|\alpha|}\int_{\delta/2<|x|< {r_4R}} |x|^{-2m+2|\alpha|-n}|D^{\alpha} u|^2 dx\notag\\
&\le& \tilde C\sum_{|\alpha|\le
2}\delta^{-4+2|\alpha|}\int_{\delta/3<|x|<\delta/2}|x|^{-2m+4-n}|D^{\alpha}u|^2 dx\notag\\
&&+C''\sum_{|\alpha|\le
2}(r_4R)^{-4+2|\alpha|}\int_{r_4R<|x|<2r_4R}|x|^{-2m+4-n}|D^{\alpha}u|^2
dx\notag\\
&\le& \tilde C'\delta^{-2m-n}\int_{|x|\le\delta}|u|^2
dx+C''(r_4R)^{-2m-n}\int_{|x|\le R_0R}|u|^2 dx,
\end{eqnarray}
where $\tilde C'$ and $C''$ are independent of $R_0$, $R$, and
$m$.

We then add $m^2(2\delta)^{-2m-n}\int_{|x|\le\delta/2}|u|^2dx$ to
both sides of \eqref{3.4} and obtain
\begin{eqnarray}\label{3.5}
&&\frac{1}{2}m^2(2\delta)^{-2m-n}\int_{|x|\le 2\delta}|u|^2dx+m^2(R_0^2R)^{-2m-n}\int_{|x|\le R_0^2R}|u|^2dx\notag\\
&=&\frac{1}{2}m^2(2\delta)^{-2m-n}\int_{|x|\le 2\delta}|u|^2dx+m^2(R_0^2R)^{-2m-n}\int_{|x|\le 2\delta}|u|^2dx\notag\\
&&+m^2(R_0^2R)^{-2m-n}\int_{2\delta<|x|\le R_0^2R}|u|^2dx\notag\\
&\le&\frac{1}{2}m^2(2\delta)^{-2m-n}\int_{|x|\le 2\delta}|u|^2dx+\frac{1}{2}m^2(2\delta)^{-2m-n}\int_{|x|\le 2\delta}|u|^2dx\notag\\
&&+m^2(R_0^2R)^{-2m-n}\int_{2\delta<|x|\le R_0^2R}|u|^2dx\notag\\
&\le&(\tilde
C'+m^2)\delta^{-2m-n}\int_{|x|\le\delta}|u|^2dx+C''(r_4R)^{-2m-n}\int_{|x|\le
R_0R}|u|^2 dx\notag\\
&=&(\tilde
C'+m^2)\delta^{-2m-n}\int_{|x|\le\delta}|u|^2dx\notag\\
&&+m^2(R_0^2R)^{-2m-n}C''m^{-2}(\frac{R_0^2}{r_4})^{2m+n}\int_{|x|\le
R_0R}|u|^2 dx.
\end{eqnarray}

We first observe that
\begin{eqnarray*}
&&
C''m^{-2}(\frac{R_0^2}{r_4})^{2m+n}=C''m^{-2}\left(\frac{4R_0}{R_0+1}\right)^{2m+n}\notag\\
&\le&C''m^{-2}(4R_0)^{2m+n}\ \le \exp(-2m)
\end{eqnarray*}
for all $R_0\le 1/16$ and $m^2\ge C''$. Thus, we obtain that
\begin{eqnarray}\label{3.6}
&&\frac{1}{2}m^2(2\delta)^{-2m-n}\int_{|x|\le 2\delta}|u|^2dx+m^2(R_0^2R)^{-2m-n}\int_{|x|\le R_0^2R}|u|^2dx\notag\\
&\leq&(\tilde
C'+m^2)\delta^{-2m-n}\int_{|x|\le\delta}|u|^2dx\notag\\
&&+m^2(R_0^2R)^{-2m-n}\exp(-2m)\int_{|x|\le R_0R}|u|^2 dx.
\end{eqnarray}

It should be noted that \eqref{3.6} is valid for all $m=j+\frac
12$ with $j\in{\mathbb N}$ and $j\ge j_0$, where $j_0$ depends on
$R_0$. Setting $R_j=(\gamma(j+\frac 12))^{-1}$ and using the
relation $m=(\gamma R)^{-1}$, we get from \eqref{3.6} that
\begin{eqnarray}\label{3.7}
&&\frac{1}{2}m^2(2\delta)^{-2m-n}\int_{|x|\le 2\delta}|u|^2dx+m^2(R_0^2R_j)^{-2m-n}\int_{|x|\le R_0^2R_j}|u|^2dx\notag\\
&\leq&(\tilde
C'+m^2)\delta^{-2m-n}\int_{|x|\le\delta}|u|^2dx\notag\\
&&+m^2(R_0^2R_j)^{-2m-n}\exp(-2cR_j^{-1})\int_{|x|\le R_0R_j}|u|^2
dx
\end{eqnarray}
for all $j\ge j_0$ and $c=\gamma^{-1}$. We now observe that
$$
R_{j+1}<R_j<2R_{j+1}\quad\text{for all}\quad j\in{\mathbb N}.
$$
Thus, if $R_{j+1}<R\leq R_j$, we can conclude that
\begin{eqnarray}\label{3.8}
\begin{cases}
&\int_{|x|\leq R_0^2R} |u|^2 dx
\leq \int_{|x|\leq R_0^2R_j} |u|^2 dx,\\
&\exp(-2cR_j^{-1})\int_{|x|\leq R_0R_j} |u|^2 dx
\leq \exp(-cR^{-1})\int_{|x|\leq R} |u|^2 dx,
\end{cases}
\end{eqnarray}
where we have used the inequality $R_0R_j\le 2R_{j+1}/16<R_{j+1}$
to derive the second inequality above. Namely, we have from
\eqref{3.7} and \eqref{3.8} that
\begin{eqnarray}\label{3.9}
&&\frac{1}{2}m^2(2\delta)^{-2m-n}\int_{|x|\le 2\delta}|u|^2dx
+m^2(R_0^2R_j)^{-2m-n}\int_{|x|\le R_0^2R}|u|^2dx\notag\\
&\leq&(\tilde
C'+m^2)\delta^{-2m-n}\int_{|x|\le\delta}|u|^2dx\notag\\
&&+m^2(R_0^2R_j)^{-2m-n}\exp(-cR^{-1})\int_{|x|\le R}|u|^2 dx.
\end{eqnarray}

If there exists $s\in{\mathbb N}$ such that
\begin{equation}\label{3.10}
R_{j+1}<R_0^{2s}\le R_j\quad\text{for some}\quad j\ge j_0,
\end{equation}
then replacing $R$ by $R_0^{2s}$ in \eqref{3.9} leads to
\begin{eqnarray}\label{3.11}
&&\frac{1}{2}m^2(2\delta)^{-2m-n}\int_{|x|\le 2\delta}|u|^2dx
+m^2(R_0^2R_j)^{-2m-n}\int_{|x|\le R_0^{2s+2}}|u|^2dx\notag\\
&\leq&(\tilde
C'+m^2)\delta^{-2m-n}\int_{|x|\le\delta}|u|^2dx\notag\\
&&+m^2(R_0^2R_j)^{-2m-n}\exp(-cR_0^{-2s})\int_{|x|\le
R_0^{2s}}|u|^2 dx.
\end{eqnarray}
Here $s$ and $R_0$ are yet to be determined. The trick now is to
find suitable $s$ and $R_0$ satisfying \eqref{3.10} and the
inequality
\begin{equation}\label{3.21}
\exp(-cR_0^{-2s})\int_{|x|\leq R_0^{2s}} |u|^2
dx\leq\frac{1}{2}\int_{|x|\leq R_0^{2s+2}} |u|^2 dx
\end{equation}
holds with such choices of $s$ and $R_0$.

It is time to use the three-sphere inequality \eqref{1.2}. To this
end, we choose $r_1=R_0^{2k+2}$, $r_2=R_0^{2k}$ and
$r_3=R_0^{2k-2}$ for $k\ge 1$. Note that $r_1/r_3<r_2/r_3\le
R_0^{2}\le 1/4$. Thus \eqref{1.2} implies
\begin{equation}\label{3.12}
\int_{|x|<R_0^{2k}}|u|^2dx/\int_{|x|<R_0^{2k+2}}|u|^2dx\leq
C^{1/\tau}(\int_{|x|<R_0^{2k-2}}|u|^2dx/\int_{|x|<R_0^{2k}}|u|^2dx)^{a},
\end{equation}
where
$$C=\max\{C_0R_0^{-2n},\exp(\beta_0(-1-4\log R_0))\}$$ and
\begin{eqnarray*}
a=\frac{1-\tau}{\tau}=\frac{A}{B}&=&\frac{(\log (r_1/r_3)-1)^2-(\log
(r_2/r_3))^2}{-1-2\log (r_2/r_3)}\\
&=&\frac{(4\log R_0-1)^2-(2\log R_0)^2}{-1-4\log R_0}.
\end{eqnarray*}
It is not hard to see that
\begin{equation}\label{3.13}
\begin{cases}
1<C\le C_0 R_0^{-\beta_1},\\
2< a\le -4\log R_0,
\end{cases}
\end{equation}
where $\beta_1=\max\{2n,4\beta_0\}$. Combining \eqref{3.13} and
using \eqref{3.12} recursively, we have that
\begin{eqnarray}\label{3.14}
&&\int_{|x|\leq R_0^{2s}} |u|^2 dx/\int_{|x|\leq R_0^{2s+2}} |u|^2
dx\notag\\
&\leq&C^{1/\tau}(\int_{|x|<R_0^{2s-2}}|u|^2dx/\int_{|x|<R_0^{2s}}|u|^2dx)^{a}\notag\\
&\leq&
C^{\frac{a^{s-1}-1}{\tau(a-1)}}(\int_{|x|<R_0^{2}}|u|^2dx/\int_{|x|<R_0^{4}}|u|^2dx)^{a^{s-1}}
\end{eqnarray}
for all $s\ge 1$. Now from the definition of $a$, we have
$\tau=1/(a+1)$ and thus
$$
\frac{a^{s-1}-1}{\tau(a-1)}=\frac{a+1}{a-1}(a^{s-1}-1)\le
3a^{s-1}.
$$
Then it follows from \eqref{3.14} that
\begin{eqnarray}\label{3.15}
&&\int_{|x|\leq R_0^{2s}} |u|^2 dx/\int_{|x|\leq R_0^{2s+2}} |u|^2
dx\notag\\
&\leq& C^{3(-4\log
R_0)^{s-1}}(\int_{|x|<R_0^{2}}|u|^2dx/\int_{|x|<R_0^{4}}|u|^2dx)^{a^{s-1}}\notag\\
&\leq& (C_0^3(R_0)^{-3\beta_1})^{(-4\log
R_0)^{s-1}}(\int_{|x|<R_0^{2}}|u|^2dx/\int_{|x|<R_0^{4}}|u|^2dx)^{a^{s-1}}.
\end{eqnarray}
Thus, by \eqref{3.15}, we can get that
\begin{eqnarray}\label{3.16}
&&\exp(-cR_0^{-2s})\int_{|x|\leq R_0^{2s}} |u|^2 dx\notag\\
&\leq& \exp(-cR_0^{-2s})(C_0^3(R_0)^{-3\beta_1})^{(-4\log
R_0)^{s-1}}\notag\\
&&(\int_{|x|<R_0^{2}}|u|^2dx/\int_{|x|<R_0^{4}}|u|^2dx)^{a^{s-1}}\int_{|x|\leq
R_0^{2s+2}} |u|^2 dx.\notag\\
\end{eqnarray}

Let $\mu=-\log R_0$, then if $R_0\ (\le 1/16)$ is sufficiently
small, i.e., $\mu$ is sufficiently large, we can see that
$$
2t\mu >(t-1)\log(4\mu)+\log(\log C_0^3+3\beta_1\mu)-\log(c/4)
$$
for all $t\in{\mathbb N}$. In other words, we have that for $R_0$
small
\begin{equation}\label{3.17}
(C_0^3R_0^{-3\beta_1})^{(-4\log
R_0)^{t-1}}<\exp(cR_0^{-2t}/4)<(1/2)\exp(cR_0^{-2t}/2)
\end{equation}
for all $t\in{\mathbb N}$. We now fix such $R_0$ so that
\eqref{3.17} holds. The constants $m_0(R_0)$ and $j_0(R_0)$ are
fixed as well. It is a key step in our proof that we can find a
universal constant $R_0$. After fixing $R_0$, we then define a
number $t_0$, depending on $R_0$ and $u$, as
\begin{eqnarray*}
t_0&=&\inf\{t\in{\mathbb R}: t\ge(\log
2-\log(ac)+\log\log(\int_{|x|<R_0^{2}}|u|^2dx/\int_{|x|<R_0^{4}}|u|^2dx))\\
&&\qquad\qquad\qquad\qquad\times(-2\log R_0-\log a)^{-1}\}.
\end{eqnarray*}
By \eqref{3.13}, one can easily check that $-2\log R_0-\log a>0$ for
all $R_0\le 1/16$. With the choice of $t_0$, we can see that
\begin{equation}\label{3.20}
(\int_{|x|<R_0^{2}}|u|^2dx/\int_{|x|<R_0^{4}}|u|^2dx)^{a^{t-1}}\le\exp(cR_0^{-2t}/2)
\end{equation}
for all $t\ge t_0$.

Let $s_1$ be the smallest positive integer such that $s_1\geq t_0$.
If
\begin{equation}\label{3.18}
R_0^{2s_1}\le R_{j_0}=(\gamma(j_0+1/2))^{-1},
\end{equation}
then we can find a $j_1\in{\mathbb N}$ with $j_1\ge j_0$ such that
\eqref{3.10} holds, i.e.,
$$
R_{j_1+1}<R_0^{2s_1}\le R_{j_1}.
$$
On the other hand, if
\begin{equation}\label{3.19}
R_0^{2s_1}> R_{j_0},
\end{equation}
then we pick the smallest positive integer $s_2>s_1$ such that
$R_0^{2s_2}\le R_{j_0}$ and thus we can also find a
$j_1\in{\mathbb N}$ with $j_1\ge j_0$ for which \eqref{3.10}
holds. We now define
\begin{equation*}
s=\begin{cases} s_1\quad\text{if}\quad\eqref{3.18}\quad\text{holds},\\
s_2\quad\text{if}\quad\eqref{3.19}\quad\text{holds}.
\end{cases}
\end{equation*}
It is important to note that with such $s$, \eqref{3.10} is
satisfied for some $j_1$ and \eqref{3.17}, \eqref{3.20} hold.
Therefore, we set $m_1=n+2(j_1+1/2)$ and $m=(m_1-n)/2$. Combining
\eqref{3.16}, \eqref{3.17} and \eqref{3.20} yields that
\begin{eqnarray*}
&&\exp(-cR_0^{-2s})\int_{|x|\leq R_0^{2s}} |u|^2 dx\notag\\
&\leq& \exp(-cR_0^{-2s})(C_0^3(R_0)^{-3\beta_1})^{(-3\log
R_0)^{s-1}}\notag\\
&&(\int_{|x|<R_0^{2}}|u|^2dx/\int_{|x|<R_0^{4}}|u|^2dx)^{a^{(s-1)}}\int_{|x|\leq
R_0^{2s+2}} |u|^2 dx.\notag\\
&\leq&\frac{1}{2}\int_{|x|\leq R_0^{2s+2}} |u|^2 dx
\end{eqnarray*}
which is \eqref{3.21}. Using \eqref{3.21} in \eqref{3.11}, we have
that
\begin{eqnarray}\label{3.22}
&&\frac{1}{2}m^2(2\delta)^{-2m-n}\int_{|x|\le 2\delta}|u|^2dx
+\frac{1}{2}m^2(R_0^2R_{j_1})^{-2m-n}\int_{|x|\le R_0^{2s+2}}|u|^2dx\notag\\
&\leq&(\tilde C'+m^2)\delta^{-2m-n}\int_{|x|\le\delta}|u|^2dx.
\end{eqnarray}
From \eqref{3.22}, we get that
\begin{eqnarray}\label{3.23}
\frac{(m_1-n)^2}{8\tilde
C'+2(m_1-n)^2}(R_0^2R_{j_1})^{-m_1}\int_{|x|\le
R_0^{2s+2}}|u|^2dx\leq\delta^{-m_1}\int_{|x|\le\delta}|u|^2dx
\end{eqnarray}
and
\begin{eqnarray*}
\frac{1}{2}m^2(2\delta)^{-2m-n}\int_{|x|\le
2\delta}|u|^2dx\leq(\tilde
C'+m^2)\delta^{-2m-n}\int_{|x|\le\delta}|u|^2dx
\end{eqnarray*}
which implies
\begin{eqnarray}\label{3.24}
\int_{|x|\le 2\delta}|u|^2dx\leq\frac{8\tilde
C'+2(m_1-n)^2}{(m_1-n)^2}2^{m_1}\int_{|x|\le\delta}|u|^2dx.
\end{eqnarray}

The estimates \eqref{3.23} and \eqref{3.24} are valid for all
$\delta\le R_0^{2s+2}/4$. Therefore, \eqref{1.3} holds with
$R_2=R_0$.  \eqref{1.4} holds with $R_3=R_0^{2s+2}/8$ and
$C_3=\frac{8\tilde C'+2(m_1-n)^2}{(m_1-n)^2}2^{m_1}$ and the proof
is now complete.\eproof

\end{document}